\documentclass{amsart}

\def\End{\operatorname {End}}

\newtheorem{thm}{Theorem}[section]
\newtheorem{lem}[thm]{Lemma}
\newtheorem{prop}[thm]{Proposition}

\theoremstyle{definition}

\theoremstyle{remark}

\numberwithin{equation}{section}

\begin{document}


\title{The Cremmer-Gervais solution of the Yang Baxter equation}

\author{Timothy J. Hodges} 
\address{University of Cincinnati, Cincinnati, OH 45221-0025,
U.S.A.}
\email{timothy.hodges@uc.edu}
\thanks{The author was supported in part by a grant from the National Science Foundation}

\date{July 1, 1997; revised October 23, 1997}



\begin{abstract}
    A direct proof is given of the fact that the Cremmer-Gervais $R$-matrix satisfies the (Quantum) Yang-Baxter equation
\end{abstract}

\maketitle


\section{Introduction}
Let $V$ be a vector space of rank $n$ over a field $F$.
Let $c \in \End V \otimes V$ be a linear operator. Define $c_{12}, c_{23} \in \End V \otimes V \otimes V$ by $c_{12} = c \otimes Id$, $c_{23} = Id \otimes c$. Then $c$ is said to satisfy the Yang-Baxter equation (YBE) if 
$$ c_{12} c_{23} c_{12} =c_{23} c_{12}c_{23}$$
An extremely interesting solution of this equation was found by Cremmer and Gervais in their paper \cite{CG}. In its slightly more general two parameter form, it is (up to a scalar)
$$
c(e_i \otimes e_j) = \begin{cases}
q e_j \otimes e_i &\text{if $i = j$} \\
qp^{i-j}e_j \otimes e_i + \sum_{i \le k < j} (q - q^{-1})p^{i-k} e_k \otimes e_{i+j-k}
	&\text{if $i<j$}\\
q^{-1}p^{i-j}e_j \otimes e_i + \sum_{j< k <i} (q^{-1}-q)p^{i-k}e_k \otimes e_{i+j-k}
	&\text{if $i>j$}
\end{cases}
$$	
where $\{e_1, \dots , e_n \}$ is a basis for $V$ and $q$ and $p$ are non-zero elements of $F$. Taking $q = p^{n/2}$ yields the original operator given by Cremmer and Gervais. 
 The derivation of this solution used some fairly technical calculations involving chiral vertex operators and is a litte inaccessible to the non-specialist. Here we give an elementary proof of this result along the same lines as the proof in \cite{Ka} for the standard solutions of the Yang-Baxter equation.


\section{Linear combinations of solutions of the YBE}

	Suppose $f$ and $g$ are solutions of the YBE and let $\alpha, \beta \in F$. Expanding out the equation 
$$
 (\alpha f + \beta g)_{12} (\alpha f + \beta g)_{23} (\alpha f + \beta g)_{12} =(\alpha f + \beta g)_{23} (\alpha f + \beta g)_{12}(\alpha f + \beta g)_{23}
$$
we see that $c = \alpha f + \beta g$ will be a solution of the YBE for all $\alpha, \beta \in A$ if the following two conditions are satisfied:
$$
 f_{12} g_{23} g_{12} +g_{12} f_{23} g_{12}+g_{12} g_{23} f_{12}
=f_{23} g_{12}g_{23}+g_{23} f_{12}g_{23}+g_{23} g_{12}f_{23}
$$
$$
 g_{12} f_{23} f_{12} +f_{12} g_{23} f_{12}+f_{12} f_{23} g_{12}
=g_{23} f_{12}f_{23}+f_{23} g_{12}f_{23}+f_{23} f_{12}g_{23}
$$
In the case where $f$ is the permutation  operator $P(e_i \otimes e_j) = e_j \otimes e_i$, the second condition is true for any $g$ (since $g_{12} P_{23} P_{12} =P_{23} P_{12}g_{23} $ and similar equalities hold for the other terms). Thus we obtain the following simple condition which we shall refer to as the {\em compatibility condition}.

\begin{lem} \label{cc} Suppose that $g\in \End V \otimes V$ is a solution of the YBE. Then $c=\alpha P + \beta g$ will be a solution of the YBE for all $\alpha, \beta \in F$ if
\begin{equation}\label{cce}
 g_{12} g_{23} P_{12} + g_{12} P_{23} g_{12}+ P_{12} g_{23} g_{12}
= g_{23} g_{12}P_{23} + g_{23} P_{12}g_{23}+P_{23} g_{12}g_{23}.
\end{equation}
\end{lem}

We shall apply this result to the case where
\begin{equation}
g(e_i \otimes e_j)= \sum_k \eta(i,j,k) e_k \otimes e_{i+j-k}
\end{equation}
and 
\begin{equation}\label{cgeta}
\eta(i,j,k) = \begin{cases}
	1 & \text{ if $i\leq k < j$} \\
	-1 & \text{ if $j\leq k < i$} \\
0 & \text{ otherwise }
	\end{cases}  
\end{equation}
Taking $\alpha = q $ and $\beta = (q - q^{-1})$ yields
$$
c(e_i \otimes e_j) = \begin{cases}
q e_j \otimes e_i &\text{if $i = j$} \\
qe_j \otimes e_i + \sum_{i \le k < j} (q - q^{-1}) e_k \otimes e_{i+j-k} &\text{if $i<j$}\\
q^{-1}e_j \otimes e_i + \sum_{j< k <i} (q^{-1}-q)e_k \otimes e_{i+j-k}&\text{if $i>j$}
\end{cases}
$$	
which is the Cremmer-Gervais operator described in the introduction in the case where $p =1$. Once we have shown that this operator satisfies the Yang-Baxter equation, it follows from some well-known ``twisting'' results \cite{H1} that the more general operator is also a solution.


\section{The compatibility condition}

In this section we check the compatibility condition (\ref{cc}) for the operator $g$ given above.

\begin{lem} Let $g \in \End V \otimes V$ be an operator of the form
$$
g(e_i \otimes e_j)= \sum_k  \eta(i,j,k) e_k \otimes e_{i+j-k}
$$
where $\eta(i,j,k) = 0$ if $k$ is not between $i$ and $j$. Then the condition of Lemma \ref{cc} is satisfied if and only if 
\begin{multline} \label{cond1}
 \eta(i,k,a+b-j)\eta(j,a+b-j,a) + \eta(i,j,b+a-k)\eta(b+a-k,k,a)\\
 + \eta(i,j,b)\eta(i+j-b,k,a) = \eta(i,k,a) \eta(i+k-a,j,b) \\
+ \eta(j,k,a) \eta(i,j+k-a,b) + \eta(j,k,j+k-b)\eta(i,j+k-b,a) 
\end{multline}
for all $i,j,k,a,b \in \{1,2, \dots, n \}$.
\end{lem}

\begin{proof}
Let $d_l =  g_{12} g_{23} P_{12} + g_{12} P_{23} g_{12}+ P_{12} g_{23} g_{12}$ and 
$d_r = g_{23} g_{12}P_{23} + g_{23} P_{12}g_{23}+P_{23} g_{12}g_{23}$. Denote 
$e_i \otimes e_j \otimes e_k$ by $[ijk]$. Then
\begin{align*}
d_l[ijk] &= \sum_{s,t} \eta(i,k,t)\eta(j,t,s)[s,j+t-s,i+k-t] \\
&+ \sum_{s,t} \eta(i,j,s)\eta(s,k,t)[t,s+k-t,i+j-s] \nonumber \\
&+ \sum_{s,t} \eta(i,j,s)\eta(i+j-s,k,t)[t,s,i+j+k-s-t] \nonumber
\end{align*}
and 
\begin{align*}
d_r[ijk]& = \sum_{s,t}  \eta(i,k,s)\eta((k+i-s,j,t)[s,t,i+j+k-s-t]\\
&+ \sum_{s,t}  \eta(j,k,s)\eta(i,j+k-s,t)[s,t,i+j+k-t-s] \nonumber\\
&+ \sum_{s,t}\eta(j,k,s)\eta(i,t,s)[s,j+k-t,i+t-s] \nonumber 
\end{align*}
Comparing the coefficients of $[a,b,i+j+k-a-b]$ then yields the result.
\end{proof}

Now set 
$$
u(x) = \begin{cases}
1 & \text{ if $x \geq 0$} \\
0 & \text{ if $x > 0$ }
	\end{cases}  
$$
and
$$
\delta(x) = \begin{cases}
1 & \text{ if $x = 0$} \\
0 & \text{ if $x \neq 0 $ }
	\end{cases}.  
$$
Notice that $\eta(i,j,k) = u(k-i)-u(k-j)$.

\begin{lem} \label {uid} For any integers $a,b,i,j,k$,
	\begin{multline*}
	u(a+b-i-j)(u(a-j)+u(b-i) -u(b-j)-u(j-b)) + u(k-b)u(a+b-i-k))\\
	=u(a-i)(u(k-b) - u(j-b) - u(b-j) +u(b+a-i-k)) + u(b-i)u(a-j)
	\end{multline*}
\end{lem}

\begin{proof}
First note that 
$$
	u(x) + u(-x) = 1 + \delta(x)
$$ and 
$$
	u(x+y)(u(x) + u(y)) = u(x)u(y) + u(x+y).
$$
From this it follows that
$$
	u(a+b-i-j)((u(a-j)+u(b-i)) = u(a-j)u(b-i) + u(a+b-i-j)
$$
$$
	u(a+b-i-j)(-u(b-j)-u(j-b)) = -u(a+b-i-j) - \delta(b-j) u(a-i)
$$
and
\begin{equation*}
\begin{split}
u(k-b)u(a+b-i-k)& =(1-u(b-k) +\delta(b-k))u(a+b-i-k)\\
	& = u(a+b-i-k)(1-u(b-k)) + \delta(b-k)u(a-i)\\
	& = u(a+b-i-k)u(a-i) - u(a-i)u(b-k) +\delta(b-k)u(a-i)\\
	& = u(a+b-i-k)u(a-i) + u(a-i)(u(k-b) -1).
\end{split}
\end{equation*}
Combining these equation yields the desired result.
\end{proof}

\begin{thm}\label{comp}
For $\eta$ as defined in (\ref{cgeta}), the operator $g(e_i \otimes e_j)= \sum_k  \eta(i,j,k) e_k \otimes e_{i+j-k}$ satisifies the compatibility condition (\ref{cce}).
\end{thm}

\begin{proof}
Expanding out the left hand side of (\ref{cond1}) using $\eta(i,j,k) = u(k-i)-u(k-j)$ yields
\begin{equation*}
	\begin{split}
	& \quad u(a+b-j-k)[u(a-k)-u(k-b)+u(j-b)-u(a-j)]\\
	& + u(a+b-i-j)[-u(j-b) + u(a-j) +u(b-i) -u(b-j)]\\
	& + u(a+b-i-k)[u(k-b) -u(a-k)]\\
	& + u(a-k)[u(b-j) -u(b-i)].
	\end{split}
\end{equation*}
Similarly the right hand side becomes
\begin{equation*}
	\begin{split}
	& \quad u(a+b-j-k)[u(a-k)-u(k-b)+u(j-b)-u(a-j)]\\
	& + u(a-k)[u(b-j) - u(b-i)-u(a+b-i-k )]\\
	& + u(a-i)[u(k-b) - u(b-j) -u(j-b)+u(a+b-i-k)]\\
	& + u(b-i)u(a-j).
	\end{split}
\end{equation*}
The equality of these two expressions follows from the identity in Lemma \ref{uid}.
\end{proof}


\section{The Yang-Baxter equation}

	In this section we verify that the operator $g(e_i \otimes e_j)= \sum_k  \eta(i,j,k) e_k \otimes e_{i+j-k}$ given above satisfies the Yang-Baxter equation. We begin by converting the problem into an identity for $\eta$.

\begin{lem} Let $g \in \End V \otimes V$ be an operator of the form
$$
g(e_i \otimes e_j)= \sum_k  \eta(i,j,k) e_k \otimes e_{i+j-k}
$$
where $\eta(i,j,k) = 0$ if $k$ is not between $i$ and $j$. Then $g$ saisfies the Yang-Baxter equation  if and only if 
\begin{multline}\label{cond2}
\sum_a \eta(j,k,a) \eta(i,a,c)\eta(i+a-c,j+k-a,h)\\
= \sum_s\eta(i,j,s)\eta(i+j-s,k,h+c-s) \eta(s,h+c-s,c)
\end{multline}
for all $i,j,k,c,h \in \{1,2, \dots, n \}$.
\end{lem}

\begin{proof}
	The left hand side of (\ref{cond2}) is the coefficient of $e_c \otimes e_h \otimes e_{i+j+k-c-h}$ in the expansion of 
$g_{23} g_{12}g_{23}(e_i \otimes e_j \otimes e_k)$. Similarly the right hand side is the coefficient of $e_c \otimes e_h \otimes e_{i+j+k-c-h}$ in 
$g_{12} g_{23}g_{12}(e_i \otimes e_j \otimes e_k)$.
\end{proof}

	The following identities are used in the proof of the next three results.

\begin{lem} \label{ids}
 For integers $a,b,c,d,e$,
\begin{enumerate}
\item $\eta(a+d,b+d,c+d) = \eta(a,b,c)$.
\item $\eta(a,b,c) = -\eta(b,a,c)$
\item $\eta(a,b,c)=\eta(-b,-a,-c-1)=\eta(a,b,a+b-c-1)$
\item $\eta(a,a+1,c) = \delta(a-c)$
\item $\sum_a\eta(b,c,a) = c-b$
\item $\eta(a,b,d) + \eta(b,c,d) = \eta(a,c,d)$
\item $\eta(a, b+1,c) \eta(c,a,b) = 0$
\item $\eta(a,b,c)\eta(c,b,d) = \eta(a,b,d)\eta(a,d+1,c)$
\item $\eta(a,b,c)\eta(d,c,e)=\eta(a,b,c)\eta(d,a,e)+ \eta(a,b,e)\eta(e+1,b,c)$
\end{enumerate}
\end{lem}

\begin{proof} The proofs are either trivial or routine calculations.
\end{proof}

\begin{lem}\label{prexi}
For any integers, $t,s,b,d,h$,
\begin{equation*}\begin{split}
&\sum_a\eta(t,s,a)\eta(b+a,d-a,h)
			= (s-t) \eta(b+t,d-t,h) \\
	& \qquad + (d-h-s) \eta(d-s,d-t,h)+(h-b-s+1)\eta(b+t,b+s,h)
\end{split}
\end{equation*}
\end{lem}

\begin{proof} Using the identities of Lemma \ref{ids},
\begin{equation*}
\begin{split}
&\sum_a\eta(t,s,a)\eta(b+a,d-a,h) \\
& = -\sum_a\eta(t,s,a)\eta(d-b-a,a,h-b)\\
& = -\sum(\eta(t,s,a)\eta(d-b-a,t,h-b) + \eta(t,s,h-b)\eta(h-b+1,s,a)) \\
& = (h-b-s+1) \eta(b+t,b+s,h) - \sum_a \eta(t,s,a)\eta(d-b-t,a,d-h-1)
\end{split}\end{equation*}
Now 
\begin{equation*}
\begin{split}
&-\sum_a \eta(t,s,a)\eta(d-b-t,a,d-h-1)\\
&=-\sum_a\eta(t,s,a)\eta(d-b-t,t,d-h-1) + \eta(t,s,d-h-1)\eta(d-h,s,a) \\
&=(s-t)\eta(t+b-d,-t,h-d) +(d-s-h)\eta(-s,-t,h-d)\\
&= (s-t) \eta(b+t,d-t,h) + (d-h-s) \eta(d-s,d-t,h)
\end{split}
\end{equation*}
Combining these two equations yields the assertion.
\end{proof}

\begin{lem} \label{xi} For any integers $i,j,k,c,h$,
\begin{equation*}
\begin{split}
&\sum_a \eta(j,k,a) \eta(i,a,c)\eta(i+a-c,j+k-a,h)\\
&= \eta(j,k,c)((k-c-1)\eta(i-c+k,j+k-c,h) + (j-h) \eta(j,j+k-c,h)\\
&\qquad + (h-i)\eta(i,i+k-c,h))\\
&+\eta(i,j,c)( (c-i+1)\eta(i+j-c,i+k-c,h) + (h-j)\eta(i+j-c,j,h)\\ &\qquad +(k-h)\eta(i+k-c,k,h)) 
\end{split}
\end{equation*}
\end{lem}

\begin{proof} By part 7 of Lemma \ref{ids}
\begin{equation*}
\begin{split}
&\sum_a \eta(j,k,a) \eta(i,a,c)\eta(i+a-c,j+k-a,h)\\
& = \sum_a (\eta(j,k,a) \eta(i,j,c)+\eta(j,k,c)\eta(c+1,k,a))\eta(i+a-c,j+k-a,h)
\end{split}
\end{equation*}
Using Lemma \ref{prexi} we obtain that
\begin{equation*}
\begin{split}
& \sum_a \eta(j,k,a) \eta(i,j,c)\eta(i+a-c,j+k-a,h)\\
&= \eta(i,j,c) ((k-j)\eta(i+j-c,k,h) + (j-h)\eta(j,k,h) \\
&\quad +(h+c-i-k+1)\eta(i+j-c,i+k-c,h))
\end{split}
\end{equation*}
and that
\begin{equation*}
\begin{split}
& \sum_a \eta(j,k,c) \eta(c+1,k,a)\eta(i+a-c,j+k-a,h)\\
&= \eta(j,k,c) ((k-c-1)\eta(i+1,j+k -c-1,h)+(j-h)\eta(j,j+k-c-1,h)\\
&\quad+(h+c-i-k+1)\eta(i+1,i+k-c,h))\\
&= \eta(j,k,c) ((k-c-1)\eta(i,j+k -c,h) +(j-h)\eta(j,j+k-c,h)\\
&\quad +(h+c-i-k+1)\eta(i,i+k-c,h))
\end{split}
\end{equation*}
Using these formulas and repeated application of the identity
$$\eta(a,b,h) + \eta(b,c,h) = \eta(a,c,h)$$
yields the result.
\end{proof}
 
\begin{thm} \label{mnthm}
For $\eta$ as defined in (\ref{cgeta}), the operator $g(e_i \otimes e_j)= \sum_k  \eta(i,j,k) e_k \otimes e_{i+j-k}$ satisifies the Yang-Baxter equation.
\end{thm}

\begin{proof}
Let 
$$\zeta(i,j,k,c,h) = \sum_a \eta(j,k,a) \eta(i,a,c)\eta(i+a-c,j+k-a,h)$$
(the left hand side of equation (\ref{cond2})). It is easily verified that the right hand side of equation (\ref{cond2}) is then $\zeta(i+j-k,i,j,h+c-k,i+j-h)$. Now
\begin{equation*}
\begin{split}
& \zeta(i+j-k,i,j,h+c-k,i+j-h) \\
&=(j+k-h-c-1) \eta(j,k,c)\eta(i-c+k,j+k-c,h) \\
&\quad +(h-j)\eta(h,j+k-c,k)\eta(i-c+k,j+k-c,h) \\
& \quad + (k-h)\eta(h,j+k-c,k)\eta(i-c+k,j+k-c,h)\\
&\quad + (h+c-i-j+1)\eta(i,j,c)\eta(i+j-c,i+k-c,h)\\
& \quad + (j-h) \eta(i+j-c,h,j) \eta(i+j-c,i+k-c,h)\\
& \quad +(h-i) \eta(i+j-c,h,i) \eta(i+j-c,i+k-c,h)
\end{split}\end{equation*}
We may then rearrange these terms one at a time using Proposition \ref{ids}:
\begin{equation*}
\begin{split}
&(j+k-h-c-1) \eta(j,k,c)\eta(i-c+k,j+k-c,h)\\
&\quad = (k-c-1) \eta(j,k,c)\eta(i-c+k,j+k-c,h) \\
& \qquad + (j-h) (j+k-h-c-1) \eta(j,k,c)\eta(i-c+k,j+k-c,h)
\end{split}
\end{equation*}
and
\begin{equation*}
\begin{split}
&(h-j)\eta(h,j+k-c,k)\eta(i-c+k,j+k-c,h)\\
& \quad = (h-j)\eta(i+k-c,j+k-c,j)\eta(i+k-c,j+1,h)\\
& \quad = (h-j)\eta(i+k-j,k,c)\eta(i+k-c,j,h).
\end{split}
\end{equation*}
Similarly
$$
(k-h)\eta(h,j+k-c,k)\eta(i-c+k,j+k-c,h)
 = (k-h)\eta(i,j,c)\eta(i+k-c,k,h),
$$
\begin{equation*}
\begin{split}
& (h+c-i-j+1)\eta(i,j,c)\eta(i+j-c,i+k-c,h)\\
&\quad = (c-i+1)\eta(i,j,c)\eta(i+j-c,i+k-c,h)\\
&\qquad + (h-j) \eta(i,j,c)\eta(i+j-c,i+k-c,h),
\end{split}
\end{equation*}
$$
 (j-h) \eta(i+j-c,h,j) \eta(i+j-c,i+k-c,h)= (j-h) \eta(i+k-j,i,c) \eta(i+k-c,j,h),
$$
and
$$
(h-i) \eta(i+j-c,h,i) \eta(i+j-c,i+k-c,h)
= (h-i) \eta(k,j,c)\eta(i+k-c,i,h)
$$
Adding these terms and rearranging easily yields $\zeta(i,j,k,c,h)$ as required.
\end{proof}

	Finally we make some observations about invertibility and the Hecke condition. Recall that $R$ is said to be {\em Hecke} if it satisfies the condition
$$
	(R-q)(R+q^{-1}) = 0
$$ for some $q$.

\begin{lem}\label{gandp}
	\begin{enumerate}
		\item $g^2=g$
		\item	$gP=-g$
		\item $Pg=g+P-I$
	\end{enumerate}
\end{lem}

\begin{proof}
	The first part follows from the identity 
$$
	\sum_k \eta(i,j,k)\eta(k,i+j-k,l) = \eta(i,j,l)
$$
which is a consequence of Lemma \ref{prexi}. The second and third parts follow from the identities
$ \eta(j,i,k) = -\eta(i,j,k)$ and $\eta(i,j,i+j-k) = \eta(i,j,k) + \delta(k-j) - \delta(k-i)$ respectively.
\end{proof}

\begin{prop} \label{hecke} Let $\alpha$ and $\beta$ be non-zero elements of $F$. The operator $R= \alpha P + \beta g$ is invertible if and only if $\alpha \neq \beta$. It is Hecke if and only if $\beta = \alpha - \alpha^{-1}$.
\end{prop}

\begin{proof}
	Using Lemma \ref{gandp} we find that 
$$
R^2 = \beta R + \alpha(\alpha-\beta)I
$$
and the proposition then follows immediately.
\end{proof}

\begin{thm} 
Let $F$ be a field  and let $V$ be a vector space with basis $\{e_1, \dots, e_n\}$. Let $c \in \End V \otimes V$ be the linear operator
$$
c(e_i \otimes e_j)= q p^{i-j} e_j \otimes e_i + \sum_k (q-q^{-1})p^{i-k} \eta(i,j,k) e_k \otimes e_{i+j-k}
$$
Then $c$ is an invertible solution of the Yang-Baxter equation.
\end{thm}

\begin{proof}
For $p=1$ the result follows from Lemma \ref{cc}, Theorem \ref{comp}, Proposition \ref{hecke} and Theorem \ref{mnthm}. For more general $p$ we apply \cite[Theorem 3.3]{H1}.
\end{proof}

\end{document}